\documentclass[11pt]{article}
\usepackage{amssymb,amsmath,latexsym}

\newtheorem{theorem}[equation]{Theorem} 
\newtheorem{corollary}[equation]{Corollary}
\newtheorem{lemma}[equation]{Lemma}
\newtheorem{proposition}[equation]{Proposition}
\newtheorem{remark}[equation]{Remark} 
\numberwithin{equation}{section}

\newcommand{\fm}{\mathfrak{m}}
\newcommand{\fn}{\mathfrak{n}}
\newcommand{\cI}{\mathcal{I}}

\newcommand{\cL}{\mathcal{L}}
\newcommand{\cU}{\mathcal{U}}
\newcommand{\qed}{\hfill $\Box$}
\newcommand{\Z}{\mathbb{Z}}
\newcommand{\proof}{{\bf Proof\ \ }}
\newcommand{\fg}{\mathfrak{g}}
\newcommand{\fh}{\mathfrak{h}}

\newcommand{\ot}{\otimes}
\newcommand{\ol}{\overline}

\newcommand{\ga}{\alpha}

\newcommand{\gd}{\delta}
\newcommand{\gl}{\lambda}
\newcommand{\gs}{\sigma}
\newcommand{\x}[2]{x_{{\ol{#1}}#2}}
\newcommand{\ad}{\hbox{ad\,}}

\title{\LARGE Representations of Multiloop Algebras}

\date{}

\author{Michael Lau\thanks{Funding from the Natural Sciences and Engineering Research Council of Canada is gratefully acknowledged. \newline
E-mail: mlau@uwindsor.ca} \\ {\small University of Windsor} \\ {\small Department of Mathematics and Statistics} \\ {\small Windsor, Ontario, Canada N9B 3P4}}

\begin{document}
\maketitle

\begin{small}
\noindent 
{\bf Abstract.}  We describe the finite-dimensional simple modules of all the (twisted and untwisted) multiloop algebras and classify them up to isomorphism.

\bigskip



\end{small}

\vskip.25truein
\section{Introduction}

Multiloop algebras are multivariable generalizations of the loop algebras appearing in affine Kac-Moody theory.  The study of these algebras and their extensions includes a substantial literature on (twisted and untwisted) multiloop, toroidal, and extended affine Lie algebras.  This paper describes the finite-dimensional simple modules of multiloop algebras and classifies them up to isomorphism.

Let $\fg$ be a finite-dimensional simple Lie algebra over an algebraically closed field $F$ of characteristic zero.  Suppose that $\gs_1,\ldots,\gs_N:\ \fg\rightarrow\fg$ are commuting automorphisms of finite orders $m_1,\ldots,m_N$, respectively.  For each $i$, fix a primitive $m_i$th root of unity $\xi_i\in F$.  Then $\fg$ decomposes into common eigenspaces relative to these automorphisms:
$$\fg=\bigoplus_{\ol{k}\in G}\fg_{\ol{k}},$$
where $\fg_{\ol{k}}=\{x\in\fg\ |\ \gs_ix=\xi_i^{k_i}x\}$ and ${\ol{k}}$ is the image of each $k\in\Z^N$ under the canonical map $\Z^N\rightarrow G=\Z/m_1\Z\times\cdots\times\Z/m_N\Z$.  The {\em multiloop algebra} of $\fg$, relative to these automorphisms, is the Lie algebra
$$\cL=\cL(\fg;\gs_1,\ldots,\gs_N)=\bigoplus_{k\in\Z^N}\fg_{\ol{k}}\ot Ft^k,$$
where $Ft^k$ is the span of $t^k=t_1^{k_1}\cdots t_N^{k_N}$, and multiplication is defined pointwise.  If the automorphisms $\gs_1,\ldots,\gs_N$ are all trivial, $\cL$ is called an {\em untwisted multiloop algebra}.  Otherwise, it is a {\em twisted multiloop algebra.}  

In the one variable case ({\em untwisted} and {\em twisted loop algebras}), the finite-dimensional simple modules can be deduced from ideas in the work of V.~Chari and A.~Pressley \cite{chari86,cp86,cp88}.  A complete list of the finite-dimensional simple modules also appears explicitly in \cite{rao93}.  These modules are classified up to isomorphism in \cite{cp87}, and a very recent paper \cite{ChFoSe} gives a detailed discussion of this problem in the twisted loop case.

A description of the finite-dimensional simple modules of the untwisted multiloop Lie algebras was first given by S.E.~Rao \cite{rao01}.  Subsequent work by P.~Batra \cite{batra} provided a complete (but redundant) list of the finite-dimensional simple modules when $\gs_1$ is a diagram automorphism and the other automorphisms $\gs_2,\ldots,\gs_N$ are all trivial.

In the one variable case, every twisted loop algebra $\cL(\fg;\tau\circ\gamma)$ defined by an inner automorphism $\tau$ and a diagram automorphism $\gamma$ is isomorphic to $\cL(\fg;\gamma)$ \cite[Prop 8.5]{kacLiebook}.  It was thus sufficient to consider only twists by diagram automorphisms in this setting.  Unfortunately, this is far from true when the number of variables is larger than one.  (See \cite[Rem 5.9]{GiPi}, for instance.)  It has recently been shown that the centreless core of almost every extended affine Lie algebra is a multiloop algebra \cite{ABFP2} (using \cite{ABFP1} and \cite{neher}).
Even for these multiloop algebras, any number of the automorphisms $\gs_i$ can be nontrivial, and any number of them can fail to be diagram automorphisms.

In this paper, we consider an arbitrary (twisted or untwisted) multiloop algebra $\cL$.  From any ideal $\cI$ of $\cL$, we construct a $G$-graded ideal $I=I(\cI)$ of the ring of Laurent polynomials $R=F[t_1^{\pm 1},\ldots,t_N^{\pm 1}]$.  If $\cI$ is the kernel of a finite-dimensional irreducible representation, the $\ol{0}$-component $I_{\ol{0}}$ of the ideal $I$ turns out to be a radical ideal of the $\ol{0}$-component of the ring $R$.  The resulting decomposition of $I_{\ol{0}}$ into an intersection of a finite number of maximal ideals produces an isomorphism 
$$\psi_a:\ \cL/\cI\rightarrow \fg\oplus\cdots\oplus \fg\quad\hbox{($r$ copies)}$$
whose composition with the quotient map $\pi:\ \cL\rightarrow \cL/\cI$ is evaluation at an $r$-tuple $a=(a_1,\ldots,a_r)$ of points $a_i\in(F^\times)^N$:
$$\psi_a\circ\pi:\ x\ot f(t)\mapsto \big(f(a_1)x,\ldots,f(a_r)x\big),$$
for any $x\ot f(t)\in\cL$.  
Since the finite-dimensional simple modules of the semisimple Lie algebra $\fg\oplus\cdots\oplus\fg$ are the tensor products of finite-dimensional simple modules for $\fg$, we obtain a complete (but redundant) list of the finite-dimensional irreducible representations of $\cL$ (Corollary \ref{main}).  Namely, any finite-dimensional simple module for $\cL$ is of the form 
$$V(\gl,a)=V_{\gl_1}(a_1)\ot \cdots \ot V_{\gl_r}(a_r),$$
where $V_{\gl_i}$ is the $\fg$-module of dominant integral highest weight $\gl_i$, and $V_{\gl_i}(a_i)$ is the $\cL$-module obtained by evaluating elements of $\cL$ at the point $a_i$, and then letting the resulting element of $\fg$ act on $V_{\gl_i}$.  The $r$-tuples $a=(a_1,\ldots,a_r)$ that occur in this process must satisfy the condition that the points $m(a_i)$ are all distinct, where $m(a_i)=(a_{i1}^{m_1},\ldots,a_{iN}^{m_N})$ is determined by the orders $m_1,\ldots,m_N$ of the automorphisms $\gs_1,\ldots,\gs_N$.  Conversely, the $\cL$-module $V(\gl,a)$ is finite-dimensional and simple if the $a_i$ satisfy this condition (Theorem \ref{evsaresimple}).

In the second half of the paper, we establish necessary and sufficient conditions for $\cL$-modules $V(\gl,a)$ and $V(\mu,b)$ to be isomorphic.  Namely, we ``pull back'' a triangular decomposition $N_-\oplus H\oplus N_+$ of $\fg\oplus\cdots\oplus\fg$ to a triangular decomposition $\psi_a^{-1}(N_-)\oplus \psi_a^{-1}(H)\oplus\psi_a^{-1}(N_+)$ of $\cL/\cI$.  The modules $V(\gl,a)$ and $V(\mu,b)$ are highest weight with respect to this decomposition of $\cL/\cI$, and they are isomorphic if and only if they have the same highest weights.  The paper concludes with three equivalent criteria for isomorphism in terms of an explicit formula (Theorem \ref{isothm}), orbits under a group action (Corollary \ref{orbitcorollary}), and equivariant maps (Corollary \ref{equivariantmapcorollary}).  These are the first such isomorphism results for modules in any multiloop setting.


One of the most interesting features is that the triangular decomposition $N_-\oplus H\oplus N_+$ is replaced with a new triangular decomposition $\psi_b\psi_a^{-1}(N_-)\oplus\psi_b\psi_a^{-1}(H)\oplus\psi_b\psi_a^{-1}(N_+)$ of $\fg^{\oplus r}$ in the computation of the highest weight of $V(\mu,b)$. Unlike diagram automorphisms, arbitrary finite order automorphisms $\gs_i$ often fail to stabilize {\em any} triangular decomposition of a finite-dimensional semisimple Lie algebra.  This fact is reflected in the change of triangular decomposition on $\fg^{\oplus r}$, and it is one of the reasons that past work considered only twists by diagram automorphisms.

The other main novelty in this classification is the passage from twists by a single nontrivial automorphism $\gs_1$ to a family of nontrivial automorphisms $\gs_1,\ldots,\gs_N$.  Here the major obstacle to past approaches was reliance on the representation theory of the fixed point subalgebra $\fg_{\ol{0}}$ under the action of the automorphisms.  While this has been a great success for the representation theory of twisted affine Lie algebras (1 variable), it is problematic when passing to a multivariable setting, since the algebra $\fg_{\ol{0}}$ is then often $0$.  We avoid this pitfall by defining the ideal $I$ in such a way (\ref{idealdef}) that whenever $\cI$ is the annihilator of a finite-dimensional module, the graded component $I_{\ol{0}}$ of $I$ is always nonzero, even if the graded component  $\fg_{\ol{0}}\ot I_{\ol{0}}$ of $\cI$ {\em is } zero.




\bigskip

\noindent
{\bf Acknowledgements:} The author thanks E.~Neher and P.~Senesi for their careful reading and corrections to an earlier draft.

\bigskip

\noindent
{\bf Note:} Throughout this paper, $F$ will be an algebraically closed field of characteristic zero.  All Lie algebras, linear spans, and tensor products will be taken over $F$ unless otherwise explicitly indicated.  We will denote the integers by $\Z$, the nonnegative integers by $\Z_+$, and the nonzero elements of $F$ by $F^\times$.

\section{Multiloop algebras and their ideals}



The following proposition is an immediate consequence of general facts about reductive Lie algebras \cite[\S 6, no. 4]{bourbakiLie}.


\begin{proposition}\label{semisimple}
Let $L$ be a perfect Lie algebra over $F$, and let $\phi:\ L\rightarrow \hbox{\em End}\,V$ be a finite-dimensional irreducible representation.  Then $L/\ker\phi$ is a semisimple Lie algebra.
\end{proposition}
\proof  The representation $\phi$ descends to a faithful representation of $L/\ker\phi$.  By \cite[Prop 6.4.5]{bourbakiLie}, any Lie algebra with a faithful finite-dimensional irreducible representation is reductive.  Moreover, $L$ is perfect.  Therefore, $L/\ker\phi$ is perfect and reductive, and hence semisimple.\qed



\bigskip
\bigskip

We now focus our attention on multiloop algebras.  Let $\fg$ be a finite-dimensional simple Lie algebra over $F$, and let $R=F[t_1^{\pm 1},\ldots,t_N^{\pm 1}]$ be the commutative algebra of Laurent polynomials in $N$ variables.  The {\em untwisted multiloop algebra} is the Lie algebra $\fg\ot R$ with (bilinear) pointwise multiplication given by 
$$[x\ot f,y\ot g]=[x,y]\ot fg$$
for all $x,y\in\fg$ and $f,g\in R$.  Suppose that $\fg$ is equipped with $N$ commuting automorphisms
$$\gs_1,\ldots,\gs_N:\ \fg\rightarrow\fg$$
of finite orders $m_1,\ldots,m_N$, respectively.  For each $i$, fix $\xi_i\in F$ to be a primitive $m_i$th root of $1$.  Then $\fg$ has a common eigenspace decomposition $\fg=\bigoplus_{\ol{k}\in G}\fg_{\ol{k}}$ where $\ol{k}$ is the image of $k=(k_1,\ldots,k_N)\in\Z^N$ under the canonical map
\begin{equation*}
\Z^N\rightarrow G=\Z/m_1\Z\times\cdots\times\Z/m_N\Z
\end{equation*}
and
\begin{equation*}
\fg_{\ol{k}}=\{x\in\fg\ |\ \gs_ix=\xi_i^{k_i}x\ \hbox{for}\ i=1,\ldots,N\}.
\end{equation*}

The {\em (twisted) multiloop algebra} $\cL=\cL(\fg;\gs_1,\ldots,\gs_N)$ is the Lie subalgebra
\begin{equation}
\cL=\bigoplus_{k\in\Z^N}\fg_{\ol{k}}\ot Ft^k\subseteq\fg\ot R,
\end{equation}
where $t^k=t_1^{k_1}\cdots t_N^{k_N}$ is multi-index notation.

Note that $R$ has a $G$-grading,
\begin{equation}\label{gradingonR}
R=\bigoplus_{\ol{k}\in G}R_{\ol{k}},
\end{equation}
where $R_{\ol 0}=F[t_1^{\pm m_1},\ldots,t_N^{\pm m_N}]$ and $R_{\ol k}=t^{k}R_{\ol 0}$ for every $k\in\Z^N$.  In this notation,
\begin{equation}\label{gradingonL}
\cL=\bigoplus_{\ol{k}\in G}\left(\fg_{\ol k}\ot R_{\ol{k}}\right).
\end{equation}

Fix an $F$-basis 
\begin{equation}\label{homobasis}
\{\x{k}{j}\ |\ j=1,\ldots,\dim\fg_{\ol{k}}\}
\end{equation}
 of $\fg_{\ol{k}}$ for all $\ol{k}\in G$.  Then 
\begin{equation}\label{deco}
\cL=\bigoplus_{\ol{k}\in G}\bigoplus_{j=1}^{\dim\fg_{\ol k}}\left( F\x{k}{j}\ot R_{\ol k}\right).
\end{equation}
Since $\fg$ is simple (hence perfect) and graded,  each $\x{k}{j}$ can be expressed as a sum of brackets of homogeneous elements $y,z\in\fg$, with $\deg y+\deg z=\ol{k}$.  For each such $k\in\Z^N$ and pair $y,z$, there exist $a,b\in\Z^N$ with $\deg y=\ol{a}$, $\deg z=\ol{b}$, and $a+b=k$.  Then the sum of the brackets $[y\ot t^a,z\ot t^b]$ will be $\x{k}{j}\ot t^{k}$.  Since these elements span $\cL$, it is clear that $\cL$ is perfect. 

Let $\pi_{\ol{k} j}$ be the projection $\pi_{\ol{k}j}:\ \cL\rightarrow F\x{k}{j}\ot R_{\ol k}$, relative to the decomposition (\ref{deco}).  We will view $\pi_{\ol{k}j}$ as a projection $\cL\rightarrow R_{\ol{k}}$ by identifying $\x{k}{j}\ot f$ with $f$ for all $f\in R_{\ol{k}}$.
Let $\cI$ be an ideal of the Lie algebra $\cL$, and let $I=I(\cI)$ be the ideal of $R$ generated by
\begin{equation}\label{idealdef}
\bigcup_{\ol{k}\in G}\bigcup_{j=1}^{\dim \fg_{\ol{k}}}\pi_{\ol{k}j}(\cI).
\end{equation}

Note that the definition of $I$ is independent of the choice of homogeneous basis $\{\x{k}{j}\}$ of $\fg$, and the ideal $I$ is $G$-graded since its generators are homogeneous with respect to the $G$-grading of $R$.  That is,
$$I=\bigoplus_{\ol{k}\in G}I_{\ol{k}}\ \ \hbox{where}\ \ I_{\ol{k}}=I\cap R_{\ol{k}}.$$
Moreover, 
$$t^{\ell-k}I_{\ol{k}}\subseteq I\cap R_{\ol{\ell}}=I_{\ol{\ell}}=t^{\ell-k}\left(t^{k-\ell}I_{\ol{\ell}}\right)\subseteq t^{\ell-k}I_{\ol{k}},$$ so
\begin{equation}\label{finegrading}
I_{\ol{\ell}}=t^{\ell-k}I_{\ol{k}},
\end{equation}
for all $k,\ell\in\Z^N$.  We will use the following technical lemma to show that $\cI=\cL\cap(\fg\ot I)$.

\begin{lemma}\label{technical} Let $\displaystyle{Y=\sum_{\ol{r}\in G}\sum_{n=1}^{\dim\fg_{\ol{r}}}\x{r}{n}\ot \pi_{\ol{r}n}(Y)\in\cI}$.  Then
$$\x{k}{i}\ot t^{k-\ell}\pi_{\ol{\ell}j}(Y)\in\cI$$
for all $k,\ell\in\Z^N$, $1\leq i\leq \dim\fg_{\ol{k}}$, and $1\leq j\leq \dim\fg_{\ol{\ell}}$.
\end{lemma}
\proof The finite-dimensional simple Lie algebra $\fg$ is a finite-dimensional simple $\fg$-module (and hence a finite-dimensional simple $\cU(\fg)$-module) under the adjoint action of $\fg$.  Fix $k,\ell\in\Z^N$, $i\in\{1,\ldots,\dim\fg_{\ol{k}}\}$, and $j\in\{1,\ldots,\dim\fg_{\ol{\ell}}\}$.  By the Jacobson Density Theorem, there exists $u\in\cU(\fg)$ such that
$$u.\x{r}{n}=\left\{
\begin{array}{ll}\x{k}{i}\ &\hbox{if}\ \ol{r}=\ol{\ell}\ \hbox{and}\ n=j\\
0\ &\hbox{otherwise,}
\end{array}\right.$$
for all $\ol{r}\in G$ and $n\in\{1,\ldots,\fg_{\ol{r}}\}$.
By the Poincar\'e-Birkhoff-Witt Theorem, write $u=\sum_{s=1}^ap_s$, where $a\geq 1$ and each $p_s$ is a monomial in the variables in $\{\x{r}{n}\ |\ \ol{r}\in G,\ n=1,\ldots,\dim\fg_{\ol{r}}\}$.  Considering the induced $G$-grading of $\cU(\fg)$, we can assume that each $p_s$ is homogeneous of degree $\ol{k-\ell}$.  Write
$$p_s=c_s\prod_{\ol{r}\in G}\prod_{n=1}^{\dim\fg_{\ol{r}}}\left(\x{r}{n}\right)^{b_{\ol{r}n}^{(s)}}$$ 
where $c_s\in F$ and $b_{\ol{r}n}^{(s)}\in\Z_+$.  

Since $p_s$ is homogeneous of degree $\ol{k-\ell}$ in the $G$-grading of $\cU(\fg)$, we can choose $\ga(s,\ol{r},n,1),\ga(s,\ol{r},n,2),\ldots,\ga(s,\ol{r},n,b_{\ol{r}n}^{(s)})\in\Z^N$ for each $s\in\{1,\ldots,a\}$, $\ol{r}\in G$, and $n\in\{1,\ldots,\dim\fg_{\ol{r}}\}$, such that
\begin{enumerate}
\item[{\rm(i)}]$\ol{r}=\ol{\ga(s,\ol{r},n,1)}=\cdots=\ol{\ga(s,\ol{r},n,b_{\ol{r}n}^{(s)})}$,
\item[{\rm (ii)}]$\displaystyle{\sum_{\ol{r}\in G}\sum_{n=1}^{\dim\fg_{\ol{r}}}\sum_{b=1}^{b_{\ol{r}n}^{(s)}}\ga(s,\ol{r},n,b)=k-\ell.}$
\end{enumerate}
Then 
$$\widetilde{p}_s=c_s\prod_{\ol{r}\in G}\prod_{n=1}^{\dim\fg_{\ol{r}}}\prod_{b=1}^{b_{\ol{r}n}^{(s)}}\left(\x{r}{n}\ot t^{\ga(s,\ol{r},n,b)}\right)$$
is in the universal enveloping algebra $\cU(\cL)$ of $\cL$, which acts on $\cI$ via the adjoint action of $\cL$ on $\cI$, and 
$$\sum_{s=1}^a\widetilde{p}_s.Y=\x{k}{i}\ot t^{k-\ell}\pi_{\ol{\ell}j}(Y).$$
Thus $\x{k}{i}\ot t^{k-\ell}\pi_{\ol{\ell}j}(Y)\in\cI$.\qed

\bigskip

\begin{proposition}\label{idealsareintersections}In the notation introduced above,
\begin{eqnarray}
\cI&=&\cL\cap(\fg\ot I)\label{intersection}\\
&=&\bigoplus_{\ol{k}\in G}\fg_{\ol{k}}\ot I_{\ol{k}}.\label{gradedideal}
\end{eqnarray}
\end{proposition}

\proof The second equality (\ref{gradedideal}) and the inclusion $\cI\subseteq\cL\cap(\fg\ot I)$ are clear, so it remains only to verify the reverse inclusion
$$\cL\cap(\fg\ot I)\subseteq\cI.$$
In light of (\ref{gradedideal}), it suffices to show that $\x{k}{i}\ot f\in\cI$ for all $\ol{k}\in G$, $i\in\{1,\ldots,\dim\fg_{\ol{k}}\}$, and $f\in I_{\ol{k}}$.

By the definition of $I$, there exist $Y_{\ol{\ell}j}\in\cI$ and $f_{\ol{\ell}j}\in R_{\ol{k-\ell}}$ such that $f=\sum_{\ol{\ell}\in G}\sum_{j=1}^{\dim\fg_{\ol{\ell}}}f_{\ol{\ell}j}\pi_{\ol{\ell}j}\left(Y_{\ol{\ell}j}\right).$  By Lemma \ref{technical},
$$\x{k}{i}\ot t^r\pi_{\ol{\ell}j}\left(Y_{\ol{\ell}j}\right)\in\cI$$
for all $r,\ell\in\Z^N$ satisfying $\ol{r}=\ol{k-\ell}$.  Since each $f_{\ol{\ell}j}\in R_{\ol{k-\ell}}$ is an $F$-linear combination of $\{t^r\ |\ \ol{r}=\ol{k-\ell}\}$, we see that 
$$\x{k}{i}\ot f_{\ol{\ell}j}\pi_{\ol{\ell}j}\left(Y_{\ol{\ell}j}\right)\in\cI$$
for all $\ol{\ell}\in G$ and $j=1,\ldots,\dim\fg_{\ol{\ell}}$.  Thus $\x{k}{i}\ot f\in\cI$.\qed

\bigskip
\bigskip

We close this section by considering the structure of $I_{\ol{0}}\subseteq R_{\ol{0}}$ in the case where $\cI$ is the kernel of an irreducible finite-dimensional representation of $\cL$.  Clearly $I_{\ol{0}}$ is an ideal of $R_{\ol{0}}$.  Moreover, it is a radical ideal:

\begin{proposition}\label{radical}
Let $\phi:\ \cL\rightarrow\hbox{\em End}\,V$ be a finite-dimensional irreducible representation of the  multiloop algebra $\cL$, and let $\cI=\ker \phi$.  Define $I=I(\cI)\subseteq R$ as above.  Then the graded component $I_{\ol{0}}$ is a radical ideal of $R_{\ol{0}}$. 
\end{proposition}
\proof  Suppose $p$ is an element of $\sqrt{I_{\ol{0}}}$, the radical of the ideal $I_{\ol{0}}=I\cap R_{\ol{0}}$ of $R_{\ol{0}}$.  Choose $k\in\Z^N$ so that $\fg_{\ol{k}}\neq 0$, and let $x\in\fg_{\ol{k}}$ be a nonzero element.

For $y\ot f\in\cL$, let $\langle y\ot f\rangle\subseteq\cL$ be the ideal (of $\cL$) generated by $y\ot f$.  Let $J=\langle x\ot t^kp\rangle$, and note that the $n$th term $J^{(n)}$ in the derived series of $J$ satisfies
$$J^{(n)}\subseteq \cL\cap\left(\fg\ot\langle p^n\rangle\right)$$
where $\langle p^n\rangle$ is the principal ideal of $R$ generated by $p^n$.  Since $I_{\ol{\ell}}=t^{\ell}I_{\ol{0}}$ for all $\ell\in\Z^N$ by (\ref{finegrading}), and since $p^n\in I_{\ol{0}}$ for $n$ sufficiently large, we see that
$$J^{(n)}\subseteq\cL\cap \left(\fg\ot I\right)$$
for $n\gg 0$.  Then by Proposition \ref{idealsareintersections}, $J^{(n)}\subseteq\cI$, so
$$\frac{J+\cI}{\cI}\subseteq\hbox{Rad}\left(\cL/\cI\right).$$

Since $\hbox{Rad}\left(\cL/\cI\right)=0$ by Proposition \ref{semisimple}, we see that $x\ot t^kp\in\cI$.  That is, $p=t^{-k}(t^kp)\in t^{-k}I_{\ol{k}}=I_{\ol{0}},$ and thus $\sqrt{I_{\ol{0}}}=I_{\ol{0}}$.\qed

\section{Some commutative algebra}

In this short section, we recall some basic commutative algebra that will be useful in the context of classifying modules for multiloop algebras.  Recall that $F$ is an algebraically closed field of characteristic zero, $F^\times=F\setminus 0$ is its group of units, and $R=F[t_1^{\pm 1},\ldots ,t_N^{\pm 1}]$.  For any ideal $I\subseteq R$, let $\mathcal{V}(I)=\left\{x\in (F^{\times})^N\ |\ f(x)=0\ \hbox{for all}\ f\in I\right\}$ be the (quasiaffine) variety corresponding to $I$, and let $\hbox{Poly}(S)=\{g\in R\ |\ g(s)=0\ \hbox{for all}\ s\in S\}$ be the ideal associated with any subset $S\subseteq (F^\times)^N$.

\begin{proposition}\label{nullstellensatz}
Let $I$ be an ideal of $R=F[t_1^{\pm 1},\ldots,t_N^{\pm 1}]$.  Then
$$\hbox{\em Poly}(\mathcal{V}(I))=\sqrt{I}.$$
\end{proposition}
\proof It is straightforward to verify that the usual proofs of the Hilbert Nullstellensatz (cf. \cite[p.~85 ]{AtiMac}, for instance) also hold for this Laurent polynomial analogue.\qed

\bigskip

The following crucial lemma is an easy consequence of the Nullstellensatz (Proposition \ref{nullstellensatz}):

\begin{lemma}
Let $J$ be a radical ideal of $R$, for which the quotient $R/J$ is a finite-dimensional vector space over $F$.  Then there exist distinct points $a_1,\ldots,a_r\in(F^\times)^N$ so that 
$$J=\fm_{a_1}\cap\cdots\cap \fm_{a_r},$$
where $\fm_{a_i}=\langle t_1-a_{i1},\ldots,t_N-a_{iN}\rangle$ is the maximal ideal corresponding to $a_i=(a_{i1},\ldots,a_{iN})$ for $i=1,\ldots ,r$.  Moreover, $\{a_1,\ldots,a_r\}$ is unique (up to permutation).
\end{lemma}
\proof  Clearly, $a\in \mathcal{V}(J)$ implies that $J\subseteq \fm_a$, so $J\subseteq \bigcap_{a\in \mathcal{V}(J)}\fm_a$.  Conversely, if $f\in\bigcap_{a\in \mathcal{V}(J)}\fm_a$ and $x\in \mathcal{V}(J)$, then $f(x)=0$ and $f\in\hbox{Poly}(\mathcal{V}(J))=\sqrt{J}=J$.  Hence $J=\bigcap_{a\in \mathcal{V}(J)}\fm_a$.


Since $J\subseteq \fm_{a_1}\cap\cdots\cap\fm_{a_r}$ for all subsets $\{a_1,\ldots,a_r\}\subseteq \mathcal{V}(J)$, we see that the ($F$-vector space) dimension of $R/(\fm_{a_1}\cap\cdots\cap\fm_{a_r})$ is bounded by $\dim_F(R/J)$.  Take a finite collection $\{a_1,\ldots ,a_r\}$ of points in $\mathcal{V}(J)$ for which this dimension is maximal.  Then $\fm_{a_1}\cap\cdots\cap\fm_{a_r}\cap\fm_{a_{r+1}}=\fm_{a_1}\cap\cdots\cap\fm_{a_r}$ for all points $a_{r+1}\in \mathcal{V}(J)$, so 
\begin{eqnarray*}
J&=&\bigcap_{b\in \mathcal{V}(J)}\fm_b\\
&=&\fm_{a_1}\cap\cdots\cap\fm_{a_r}\cap\left(\bigcap_{b\in \mathcal{V}(J)}\fm_b\right)\\
&=&\fm_{a_1}\cap\cdots\cap\fm_{a_r}.
\end{eqnarray*}

To see that $\{a_1,\ldots,a_r\}\subseteq(F^\times)^N$ is uniquely determined, suppose that $J=\fm_{a_1}\cap\cdots\cap\fm_{a_r}=\fm_{b_1}\cap\cdots\cap\fm_{b_s}$ for some $a_1,\ldots,a_r,b_1,\ldots,b_s\in (F^\times)^N$.  Then
\begin{eqnarray*}
\{a_1,\ldots ,a_r\}&=&\mathcal{V}(\fm_{a_1}\cap\cdots\cap \fm_{a_r})\\
&=&\mathcal{V}(J)\\
&=&\mathcal{V}(\fm_{b_1}\cap\cdots\cap \fm_{b_s})\\
&=&\{b_1,\ldots,b_s\}.
\end{eqnarray*}
\qed

\bigskip

Note that the ideal $I_{\ol{0}}\subseteq R_{\ol{0}}$ of Proposition \ref{radical} is radical and cofinite.  Viewing $R_{\ol{0}}=F[t_1^{\pm m_1},\ldots,t_N^{\pm m_N}]$ as the ring of Laurent polynomials in the variables $t_1^{m_1},\ldots ,t_N^{m_N}$, we see that 
\begin{equation}
I_{\ol{0}}=M_{a_1}\cap\cdots\cap M_{a_r},
\end{equation}
where $\{a_1,\ldots,a_r\}=\mathcal{V}(I_{\ol{0}})$ is a set of distinct points in $(F^\times)^N$, and $M_{a_i}=\langle t_1^{m_1}-a_{i1},\ldots,t_N^{m_N}-a_{iN}\rangle_{R_{\ol{0}}}$ is the maximal ideal of $R_{\ol{0}}$ corresponding to the point $a_i=(a_{i1},\ldots,a_{iN})$.  Then by the Chinese Remainder Theorem, we have the following corollary:

\begin{corollary}\label{chineseremainder}
Let $I_{\ol{0}}$ and $R_{\ol{0}}$ be as in Proposition \ref{radical}.  Then there exist unique (up to reordering) points $a_1,\ldots,a_r\in (F^\times)^N$ so that the canonical map 
\begin{eqnarray}
R_{\ol{0}}/I_{\ol{0}}&\rightarrow& R_{\ol{0}}/M_{a_1}\times\cdots\times R_{\ol{0}}/M_{a_r}\\
f+I_{\ol{0}}&\mapsto&(f+M_{a_1},\ldots,f+M_{a_r})
\end{eqnarray}
is a well-defined $F$-algebra isomorphism.
\end{corollary}
\qed

\section{Classification of simple modules}

We now return to classifying the finite-dimensional simple modules of  multiloop algebras.  As in \S2, let $\fg$ be a finite-dimensional simple Lie algebra, and let $\phi:\ \cL\rightarrow \hbox{End}\,V$ be a finite-dimensional irreducible representation of a  multiloop algebra $\cL=\cL(\fg;\gs_1,\ldots,\gs_N)$ defined by commuting automorphisms $\gs_1,\ldots,\gs_N:\ \fg\rightarrow\fg$ of order $m_1,\ldots,m_N$, respectively.

Letting $\cI=\ker\phi$, $I=I(\cI)$, $G=\Z/m_1\Z\times\cdots\times\Z/m_N\Z$, and $R=F[t_1^{\pm 1},\ldots,t_N^{\pm 1}]$ be defined as in \S2, we see that
$$\cL=\bigoplus_{\ol{k}\in G}\fg_{\ol{k}}\ot R_{\ol{k}}\hbox{\ \ and \ \ }\cI=\bigoplus_{\ol{k}\in G}\fg_{\ol{k}}\ot I_{\ol{k}},$$
by Proposition \ref{idealsareintersections}.  Since $\cI$ is a $G$-graded ideal of $\cL$, we have
\begin{eqnarray}
\cL/\cI&=&\bigoplus_{\ol{k}\in G}\left(\left(\fg_{\ol{k}}\ot R_{\ol{k}}\right)/\left(\fg_{\ol{k}}\ot I_{\ol{k}}\right)\right)\\
&=&\bigoplus_{\ol{k}\in G}\fg_{\ol{k}}\ot \left(R_{\ol{k}}/I_{\ol{k}}\right).
\end{eqnarray}

Each graded component $R_{\ol{k}}/I_{\ol{k}}$ of $R/I$ is an $R_{\ol{0}}$-module, and it is easy to check that the map
\begin{eqnarray}
\mu_k:\ R_{\ol{0}}/I_{\ol{0}}&\rightarrow& R_{\ol{k}}/I_{\ol{k}}\\
f+I_{\ol{0}}&\mapsto&t^kf+I_{\ol{k}}
\end{eqnarray}
is a well-defined $R_{\ol{0}}$-module homomorphism for each $k\in\Z^N$ and $f\in R_{\ol{0}}$.  By (\ref{gradingonR}) and (\ref{finegrading}), $R_{\ol{k}}=t^kR_{\ol{0}}$ and $t^{-k}I_{\ol{k}}=I_{\ol{0}}$, so the map $\mu_k$ is both surjective and injective.  Hence the following lemma holds:
\begin{lemma}\label{componentsiso}
Let $k\in\Z^N$.  Then the map $\mu_k:\  R_{\ol{0}}/I_{\ol{0}}\rightarrow R_{\ol{k}}/I_{\ol{k}}$ is a well-defined isomorphism of $R_{\ol{0}}$-modules.  In particular, each graded component $R_{\ol{k}}/I_{\ol{k}}$ has the same dimension (as a vector space):
$$\dim\left(R_{\ol{0}}/I_{\ol{0}}\right)=\dim\left(R_{\ol{k}}/I_{\ol{k}}\right).$$
\end{lemma}\qed

Let $a_1,\ldots ,a_r\in \left(F^\times\right)^N$ be the (unique) points defined by Corollary \ref{chineseremainder}, and let $b_i=(b_{i1},\ldots,b_{iN})$ be a point in $\left(F^\times\right)^N$ such that $b_{ij}^{m_j}=a_{ij}$ for all $1\leq i\leq r$ and $1\leq j\leq N$.  Recall that $I_{\ol{k}}=t^kI_{\ol{0}}$ for all $k\in\Z^N$, and $I_{\ol{0}}$ is contained in the ideal $M_{a_i}$ of $R_{\ol{0}}$ for $i=1,\ldots ,r$.  Therefore, the map
\begin{eqnarray}\label{epi}
\psi=\psi_b:\ \cL&\rightarrow&\fg\oplus\cdots\oplus\fg\ (r\ \ \ \hbox{copies})\\
x\ot f&\mapsto&\big(f(b_1)x,\ldots,f(b_r)x\big)
\end{eqnarray}
descends to a well-defined Lie algebra homomorphism
\begin{equation}\label{iso}
\ol{\psi}:\ \cL/\ker\phi\rightarrow\fg\oplus\cdots\oplus\fg.
\end{equation}

\begin{theorem}\label{explicitss}
The map $\ol{\psi}:\ \cL/\ker\phi\rightarrow\fg\oplus\cdots\oplus\fg$ in (\ref{iso}) is a Lie algebra isomorphism.
\end{theorem}
\proof Let $k\in\Z^N$, and let 
$$\ol{\psi}_{\ol{k}}:\ \fg_{\ol{k}}\ot \left(R_{\ol{k}}/I_{\ol{k}}\right)\rightarrow\fg_{\ol{k}}\oplus\cdots\oplus\fg_{\ol{k}}$$
be the restriction of $\ol{\psi}$ to the graded component $\fg_{\ol{k}}\ot \left(R_{\ol{k}}/I_{\ol{k}}\right)$ of $\cL/\ker\phi$.  

Note that the map $\ol{\psi}$ is injective if each $\ol{\psi}_{\ol{k}}$ is injective.  In the notation of (\ref{homobasis}), if 
$$u=\sum_{j=1}^{\dim\fg_{\ol{k}}}x_{\ol{k}j}\ot \left(t^kf_j(t)+ I_{\ol{k}}\right)$$
is in the kernel of $\ol{\psi}_{\ol{k}}$
for some collection of $f_j\in R_{\ol{0}}$, then $b_i^kf_j(b_i)=0$ for all $i$ and $j$.  Then for all $i$ and $j$, we have $f_j(b_i)=0$ and $f_j\in M_{a_i}$, where $ M_{a_i}$ is the ideal of $R_{\ol{0}}$ generated by $\{t_{\ell}^{m_{\ell}}-a_{i\ell}\ |\ \ell=1,\ldots,N\}$.  Hence $f_j\in\bigcap_{i=1}^r M_{a_i}=I_{\ol{0}}$, so $t^kf_j(t)\in t^kI_{\ol{0}}=I_{\ol{k}}$, and
\begin{eqnarray*}
\sum_{j=1}^{\dim \fg_{\ol{k}}}x_{\ol{k}j}\ot t^kf_j(t)&\in&\fg_{\ol{k}}\ot I_{\ol{k}}\\
&\subseteq&\ker\phi.
\end{eqnarray*}
Hence $u=0$ in $\cL/\ker\phi$, so $\ol{\psi}_{\ol{k}}$ (and thus $\ol{\psi}$) is injective.

By Lemma \ref{componentsiso}, $\dim\left(R_{\ol{\ell}}/I_{\ol{\ell}}\right)=\dim\left(R_{\ol{0}}/I_{\ol{0}}\right)$ for all $\ell\in\Z^N$.  Therefore, 
\begin{eqnarray*}
\dim\left(\cL/\ker\phi\right)&=&\sum_{\ol{\ell}\in G}\left(\dim \fg_{\ol{\ell}}\right)\left(\dim \left(R_{\ol{\ell}}/I_{\ol{\ell}}\right)\right)\\
&=&\dim\left(R_{\ol{0}}/I_{\ol{0}}\right)\dim\fg.
\end{eqnarray*}
Since $F$ is algebraically closed, $R_{\ol{0}}/ M_{a_i}\cong F$ for every $i$, so the ($F$-vector space) dimensions satisfy
\begin{eqnarray*}
\dim\left(R_{\ol{0}}/I_{\ol{0}}\right)&=&\dim\left(R_{\ol{0}}/ M_{a_1}\times\cdots\times R_{\ol{0}}/ M_{a_r}\right)\\
&=&r,
\end{eqnarray*}
by Corollary \ref{chineseremainder}.  Therefore, $\ol{\psi}$ is an injective homomorphism between two Lie algebras of equal dimension, so $\ol{\psi}$ is an isomorphism.\qed

\bigskip

The finite-dimensional simple modules over direct sums of copies of the Lie algebra $\fg$ are tensor products of finite-dimensional simple modules over $\fg$.  (See \cite[\S7, no.~7]{bourbaki}, for instance.)  We can thus conclude that the finite-dimensional simple modules for  multiloop algebras are pullbacks (under $\psi$) of tensor products of finite-dimensional simple modules over $\fg$.  

Fix a Cartan subalgebra $\fh\subset\fg$, a base $\Delta$ of simple roots, and weights $\gl_i\in\fh^*$ for $i=1,\ldots,r$.  Then we will write $V_{\gl_i}(b_i)$ for the simple $\fg$-module $V_{\gl_i}$ of highest weight $\gl_i$, equipped with the $\cL$-action given by
$$\big(x\ot f(t)\big).v=f(b_i)xv,$$
for all $x\ot f\in\cL$ and $v\in V_{\gl_i}$.  The tensor product of such a family of evaluation modules will be denoted
\begin{equation}
V(\gl,b)=V_{\gl_1}(b_1)\ot\cdots\ot V_{\gl_r}(b_r),
\end{equation}
and we will write $m(b_i)$ for the point $(b_{i1}^{m_1},\ldots,b_{i N}^{m_N})\in \left(F^\times\right)^N$ for $i=1,\cdots,r$.  We have now proved one of our main results:

\begin{corollary}\label{main}
Let $V$ be a finite-dimensional simple module for the  multiloop algebra $\cL$. 
Then there exist $b_1,\ldots,b_r\in \left(F^\times\right)^N$ and $\gl_1,\ldots,\gl_r$ dominant integral weights for $\fg$ such that 
$V\cong V(\gl,b),$
where $m(b_i)\neq m(b_j)$ whenever $i\neq j$.
\end{corollary}
\qed

Conversely, if the points $m(b_i)\in(F^\times)^N$ are pairwise distinct, then such a tensor product of evaluation modules is simple:

\begin{theorem}\label{evsaresimple}
Let $\gl_1,\ldots,\gl_r$ be dominant integral weights for $\fg$, and let $b_1,\ldots,b_r\in (F^\times)^N$ satisfy the property that $m(b_i)\neq m(b_j)$ whenever $i\neq j$.  Then $V(\gl,b)$ is a finite-dimensional simple $\cL$-module.
\end{theorem}
\proof Let $I_{\ol{0}}$ be the intersection $\displaystyle{\bigcap_{i=1}^r M_{a_i}}$ of the maximal ideals $ M_{a_i}$ of $R_{\ol{0}}$ corresponding to the points $a_i=m(b_i)$.  For any $k,\ell\in\Z^N$, we see that $t^{k-\ell}I_{\ol{0}}=I_{\ol{0}}$ if $\ol{k}=\ol{\ell}$ as elements of $G=\Z/m_1\Z\times\cdots\times\Z/m_N\Z$.  Thus $t^kI_{\ol{0}}=t^{\ell}I_{\ol{0}}$ if $\ol{k}=\ol{\ell}$, so we can unambiguously define $I_{\ol{k}}=t^kI_{\ol{0}}$ for any $k\in\Z^N$.

Since $a_1,\ldots,a_r$ are pairwise distinct points in $(F^\times)^N$, the proof of Theorem \ref{explicitss} (in particular, the appeal to Corollary \ref{chineseremainder}) shows that the map
\begin{eqnarray*}
\psi:\ \cL&\rightarrow&\fg\oplus\cdots\oplus\fg\ \ (r\ \hbox{copies})\\
x\ot f(t)&\mapsto&\big(f(b_1)x,\ldots,f(b_r)x\big)
\end{eqnarray*}
is surjective.  Then since each $V_{\gl_i}$ is a simple $\fg$-module, the tensor product $V_{\gl_1}\ot\cdots\ot V_{\gl_r}$ is a simple module over $\fg\oplus\cdots\oplus \fg$, and the pullback $V(\gl,b)$ is a simple $\cL$-module.\qed


\begin{remark} {\em
It is not difficult to verify that if $m(b_i)=m(b_j)$ for some $i\neq j$ for which $\gl_i$ and $\gl_j$ are both nonzero, then $V(\gl,b)$ is {\em not} simple.  However, as we do not need this fact for the classification of simple modules, we will omit its proof. 
}\end{remark}

\section{Isomorphism classes of simple modules}

By Corollary \ref{main} and Theorem \ref{evsaresimple}, the finite-dimensional simple modules of the  multiloop algebra $\cL(\fg;\gs_1,\ldots,\gs_N)$ are precisely the tensor products 
\begin{equation}\label{eqn1}
V(\gl,a)=V_{\gl_1}(a_1)\ot\cdots\ot V_{\gl_r}(a_r)
\end{equation}
for which all the $\gl_i\in\fh^*$ are dominant integral, and $m(a_i)\neq m(a_j)$ whenever $i\neq j$. 
If $\gl_i=0$ for some $i$, then $V_{\gl_i}(a_i)$ is the trivial module, and (up to isomorphism) this term can be omitted from the tensor product (\ref{eqn1}).  With the convention that empty tensor products of $\cL$-modules are the $1$-dimensional trivial module, we may assume that every $\gl_i$ is a {\em nonzero} dominant integral weight in (\ref{eqn1}).

\bigskip

To proceed further, we will need a lemma about how highest weights depend on triangular decompositions.

Let $L$ be a finite-dimensional semisimple Lie algebra with Cartan subalgebra $H$ and base of simple roots $\Delta\subset H^*$.
Recall that the group $\hbox{Aut}\,L$ of automorphisms of $L$ is (canonically) a semidirect product of the group  $\hbox{Int}\,L$ of inner automorphisms and the group $\hbox{Out}\,L$ of diagram automorphisms with respect to $(H,\Delta)$:
\begin{equation}\label{inNout}
\hbox{Aut}\,L=\hbox{Int}\,L\rtimes\hbox{Out}\,L.
\end{equation}  
See \cite[IX.4]{jacobson}, for instance.  Every automorphism $\theta$ can thus be decomposed as $\theta=\tau\circ\gamma$ with an {\em inner part} $\tau\in\hbox{Int}\,L$ and {\em outer part} $\gamma\in\hbox{Out}\,L$.



\bigskip

\begin{lemma}\label{5.2}
Let $H$ be a Cartan subalgebra of a finite-dimensional semisimple Lie algebra $L$, and let $\Delta\subset H^*$ be a base of simple roots.  Suppose that $V$ is a finite-dimensional simple $L$-module of highest weight $\gl$ with respect to $(H,\Delta)$, and $\theta\in\hbox{Aut}\,L$.  Write $\theta=\tau\circ\gamma$ for some $\tau\in\hbox{Int}\,L$ and $\gamma\in\hbox{Out}\,L$.

Then $\Delta\circ \theta^{-1}=\{\ga\circ\theta^{-1}\ |\ \ga\in\Delta\}$ is a base of simple roots for $L$, relative to the Cartan subalgebra $\theta(H)\subset L$, and $V$ has highest weight $\gl\circ\tau^{-1}$ with respect to $\big(\theta(H),\Delta\circ\theta^{-1}\big)$.
\end{lemma}
\proof Any diagram automorphism with respect to $(H,\Delta)$ will preserve $H$ and $\Delta$, so $V$ has highest weight $\gl$ with respect to $\big(\gamma(H),\Delta\circ \gamma^{-1}\big)=(H,\Delta)$.  Therefore, it is enough to prove the lemma for the case where $\theta=\tau$ is an inner automorphism.
Since inner automorphisms are products of automorphisms of the form $\exp(\ad x)$ for ad-nilpotent elements $x\in L$, we may also assume, without loss of generality, that $\tau=\exp(\ad u)$ for some ad-nilpotent element $u$.

Let $\rho:\ L\rightarrow \hbox{End}\,V$ be the homomorphism describing the action of $L$ on $V$.  Then for any $v\in V$,
\begin{align}\label{conj}
\tau(h).v&=\big(\exp(\ad u)(h)\big).v\\
&=e^{\rho(u)}\rho(h)e^{-\rho(u)}v,\label{conj2}
\end{align}
where $e^{\rho(u)}$ denotes the matrix exponential of the endomorphism $\rho(u)$.

The map $e^{\rho(u)}$ is invertible, so for any nonzero element
$$w\in V_\ga^H:=\{v\in V\ |\ h.v=\ga(h)v\ \hbox{for all}\ h\in H\},$$
we see that $e^{\rho(u)}w\neq 0$, and using (\ref{conj})--(\ref{conj2}),
\begin{align*}
\tau(h).e^{\rho(u)}w&=e^{\rho(u)}\rho(h)e^{-\rho(u)}e^{\rho(u)}w\\
&=\ga(h)e^{\rho(u)}w.
\end{align*}
That is, 
$$e^{\rho(u)}V_\ga^H\subseteq V_{\ga\circ\tau^{-1}}^{\tau(H)}:=\{v\in V\ |\ h.v=\ga\circ\tau^{-1}(h).v\ \hbox{for all}\ h\in\tau(H)\}.$$
The reverse inclusion follows similarly by considering $\tau^{-1}=\exp(-\ad u)$, so
\begin{equation}
e^{\rho(u)}V_\ga^H=V_{\ga\circ\tau^{-1}}^{\tau(H)}
\end{equation}
for all $\ga\in H^*$.  In the case where $V$ is the adjoint module $L$, we now see that $\ga$ is a root relative to $H$ if and only if $\ga\circ\tau^{-1}$ is a root relative to $\tau(H)$.  It follows easily that $\Delta\circ\tau^{-1}$ is a base of simple roots for $L$, with respect to the Cartan subalgebra $\tau(H)$.

The second part of the lemma also follows easily, since $V_{\gl\circ\tau^{-1}}^{\tau(H)}=e^{\rho(u)}V_\gl^H$ is nonzero, but 
$$V_{\gl\circ\tau^{-1}+\ga\circ\tau^{-1}}^{\tau(H)}=e^{\rho(u)}V_{\gl+\ga}^H=0$$
for all $\ga\in\Delta$.  That is, the highest weight of $V$ is $\gl\circ\tau^{-1}$, relative to $(\tau(H),\Delta\circ \tau^{-1})=(\theta(H),\Delta\circ\theta^{-1})$.\qed

\bigskip

Fix a base $\Delta$ of simple roots with respect to a Cartan subalgebra $\fh\subseteq \fg$.  The following theorem gives necessary and sufficient conditions for modules of the form $V(\gl,a)$ to be isomorphic.

\begin{theorem}\label{isothm}
Let $\gl=(\gl_1,\ldots,\gl_r)$ and $\mu=(\mu_1,\ldots,\mu_s)$ be sequences of nonzero dominant integral weights with respect to $\Delta$.  Suppose $a=(a_1,\ldots,a_r)$ and $b=(b_1,\ldots,b_s)$ are sequences of points in $\left(F^\times\right)^N$ with $m(a_i)\neq m(a_j)$ and $m(b_i)\neq m(b_j)$ whenever $i\neq j$.  

Then the finite-dimensional simple $\cL$-modules $V(\gl,a)$ and $V(\mu,b)$ are isomorphic if and only if $r=s$ and there is a permutation $\pi\in S_r$ satisfying the following two conditions for $i=1,\ldots,r$:
$$m(a_i)=m(b_{\pi(i)})\quad \hbox{and}\quad \gl_i=\mu_{\pi(i)}\circ\gamma_i,$$
where $\gamma_i$ is the outer part of the automorphism $\omega_i:\ \fg\rightarrow\fg$ defined by $\omega_i(x)=(b_{\pi(i)}^k/a_i^k)x$
for all $k\in\Z^N$ and $x\in\fg_{\ol{k}}$.
\end{theorem}
\proof Let $\phi_{\gl,a}:\ \cL\rightarrow \hbox{End}\,V(\gl,a)$ and $\phi_{\mu,b}:\ \cL\rightarrow \hbox{End}\,V(\mu,b)$ be the Lie algebra homomorphisms defining the representations $V(\gl,a)$ and $V(\mu,b)$.  By Theorem \ref{explicitss}, the kernel of $\phi_{\gl,a}$ is equal to the kernel of the evaluation map $\psi_a$, defined by
\begin{align*}
\psi_a:\quad\quad \cL&\longrightarrow \fg\oplus\cdots\oplus\fg\\
x\ot f&\mapsto\big(f(a_1)x,\ldots,f(a_r)x\big)
\end{align*}
for all $x\ot f\in\cL$.  Similarly, $\ker\phi_{\mu,b}=\ker\psi_b$.  

If the $\cL$-modules $V(\gl,a)$ and $V(\mu,b)$ are isomorphic, then $\ker\phi_{\gl,a}=\ker\phi_{\mu,b}$, so $\ker\psi_a=\ker\psi_b$.  But $\ker\psi_a=\bigoplus_{\ol{k}\in G}\fg_{\ol{k}}\ot I_{\ol{k}}$, where $I_{\ol{k}}=t^kI_{\ol{0}}$ for all $k\in\Z^N$, and 
$$I_{\ol{0}}=M_{m(a_1)}\cap\cdots\cap M_{m(a_r)},$$
where $M_{m(a_i)}=\langle t_1^{m_1}-a_{i1}^{m_1},\ldots,t_N^{m_N}-a_{iN}^{m_N}\rangle_{R_{\ol{0}}}$ is the maximal ideal of $R_{\ol{0}}=F[t_1^{\pm m_1},\ldots,t_N^{\pm m_N}]$ corresponding to the point $m(a_i)=(a_{i1}^{m_1},\ldots,a_{iN}^{m_N})$.  Since $\ker\psi_a=\ker\psi_b$, we see that (in the notation of \S 3):
\begin{align*}
\{m(a_1),\ldots,m(a_r)\}&=\mathcal{V}(M_{m(a_1)}\cap\cdots\cap M_{m(a_r)})\\
&=\mathcal{V}(I_{\ol{0}})\\
&=\mathcal{V}(M_{m(b_1)}\cap\cdots\cap M_{m(b_s)})\\
&=\{m(b_1),\ldots,m(b_s)\}.
\end{align*}
Hence $r=s$, and there is a permutation $\pi\in S_r$ such that $m(a_i)=m(b_{\pi(i)})$ for $i=1,\dots,r$.  We will write $\pi(b)=(b_{\pi(1)},\ldots,b_{\pi(r)})$.

Let 
$\fg=\fn_-\oplus\fh\oplus\fn_+$ be the triangular decomposition of $\fg$ relative to $\Delta$.  Assuming that $r=s$ and $m(a_i)=m(b_{\pi(i)})$ for all $i$, view $V_\gl=V_{\gl_1}\ot \cdots\ot V_{\gl_r}$ and $V_{\pi(\mu)}=V_{\mu_{\pi(1)}}\ot\cdots\ot V_{\mu_{\pi(r)}}$ as highest weight modules of the semisimple Lie algebra $\fg^{\oplus r}$, relative to the triangular decomposition
\begin{equation}\label{g-tridecomp}
\fg^{\oplus r}=\big(\fn_-^{\oplus r}\big)\oplus \big(\fh^{\oplus r}\big) \oplus \big(\fn_+^{\oplus r}\big).
\end{equation}
The highest weights of $V_\gl$ and $V_{\pi(\mu)}$ are $\gl$ and $\pi(\mu)=(\mu_{\pi(1)},\ldots,\mu_{\pi(r)})$, respectively, where
$\gl(h_1,\ldots,h_r)=\sum_i\gl_i(h_i)$ for all $(h_1,\ldots, h_r)\in\fh^{\oplus r}$, and $\pi(\mu)\in\big(\fh^{\oplus r}\big)^*$ is defined analogously.

We can pull back the triangular decomposition (\ref{g-tridecomp}) via the isomorphism $\ol{\psi}_a:\ \cL/\ker\psi_a\rightarrow\fg^{\oplus r}$ defined in (\ref{iso}).  Then $V(\gl,a)$ and $V(\mu,b)$ are irreducible highest weight modules of the semisimple Lie algebra $\cL/\ker\psi_a$, relative to the triangular decomposition
\begin{equation}\label{L-tridecomp}
\cL/\ker\psi_a=\ol{\psi}_a^{-1}\big(\fn_-^{\oplus r}\big)\oplus\ol{\psi}_a^{-1}\big(\fh^{\oplus r}\big)\oplus\ol{\psi}_a^{-1}\big(\fn_+^{\oplus r}\big).
\end{equation}
The $\cL$-modules $V(\gl,a)$ and $V(\mu,b)$ are isomorphic if and only if they have the same highest weights relative to the decomposition (\ref{L-tridecomp}).  Since $\ol{\psi}_a$ maps the decomposition (\ref{L-tridecomp}) to the decomposition (\ref{g-tridecomp}), the highest weight of $V(\gl,a)$ is clearly $\gl\circ\ol{\psi}_a:\ \ol{\psi}_a^{-1}\big(\fh^{\oplus r}\big)\rightarrow F$.  

The highest weight of $V(\mu,b)$ is $\nu\circ\ol{\psi}_{\pi(b)}$, where $\nu\in\big(\ol{\psi}_{\pi(b)}\ol{\psi}_a^{-1}(\fh^{\oplus r})\big)^*$ is the highest weight of $V_{\pi(\mu)}$ relative to the new triangular decomposition
\begin{equation}\label{g-tridecomp2}
\fg^{\oplus r}=\ol{\psi}_{\pi(b)}\ol{\psi}_a^{-1}\big(\fn_-^{\oplus r}\big)\oplus\ol{\psi}_{\pi(b)}\ol{\psi}_a^{-1}\big(\fh^{\oplus r}\big)\oplus\ol{\psi}_{\pi(b)}\ol{\psi}_a^{-1}\big(\fn_+^{\oplus r}\big).
\end{equation}
Let $\ol{\psi}_{\pi(b)}\ol{\psi}_a^{-1}=\tau\circ\gamma$ be a decomposition into an inner automorphism $\tau$ and a diagram automorphism $\gamma$ with respect to $\big(\fh^{\oplus r},\Delta\big)$.  By Lemma \ref{5.2}, $\nu=\pi(\mu)\circ\tau^{-1}$, so the two modules $V(\gl,a)$ and $V(\mu,b)$ are isomorphic if and only if $\gl\circ\ol{\psi}_a=\pi(\mu)\circ\tau^{-1}\circ\ol{\psi}_{\pi(b)}$ on $\ol{\psi}_a^{-1}(\fh^{\oplus r})$.  That is, $V(\gl,a)\cong V(\mu,b)$ if and only if
\begin{equation}
\gl=\pi(\mu)\circ\gamma
\end{equation}
on $\fh^{\oplus r}$.  To finish the proof, it is enough to write down an explicit formula for the automorphism $\ol{\psi}_{\pi(b)}\ol{\psi}_a^{-1}=\tau\circ\gamma$ of $\fg^{\oplus r}$.

For each $x\in\fg$, let 
$x^i=(0,\ldots,x,\ldots,0)\in\fg^{\oplus r}\quad \hbox{($x$ in the $i$th position).}$
If $k\in\Z^N$ and $x\in\fg_{\ol{k}}$, then we see that 
\begin{align*}
\ol{\psi}_a^{-1}(x^i)&=a_i^{-k}x\ot t^kf_i(t)+\ker\psi_a
\end{align*}
in $\cL/\ker\psi_a$, for any $f_i(t)\in R_{\ol{0}}$ with $f_i(a_j)=\gd_{ij}$ for all $j=1,\ldots,r$.  Since $f_i\in R_{\ol{0}}=F[t_1^{\pm m_1},\ldots,t_N^{\pm m_N}]$ and $m(a_j)=m(b_{\pi(j)})$ for all $j$, we see that $f_i(b_{\pi(j)})=\gd_{ij}$, and
\begin{align*}
\ol{\psi}_{\pi(b)}\ol{\psi}_a^{-1}(x^i)&=\left(\frac{b_{\pi(i)}^k}{a_i^k}\right)x^i.
\end{align*}
\qed

\bigskip

Theorem \ref{isothm} may also be interpreted in terms of a group action on the space of parameters $(\gl,a)$ defining the finite-dimensional simple modules of $\cL$.  Let $G^r=G\times\cdots\times G$ ($r$ factors), where $G$ is the finite abelian group $G=\langle \gs_1\rangle\times\cdots\times\langle \gs_N\rangle$, as before.  Note that $G$ acts on $\big(F^\times\big)^N$ via the primitive $m_i$th roots of unity $\xi_i$ used in the definition of $\cL$:
$$(\gs_1^{c_1},\ldots,\gs_N^{c_N}).(d_1,\ldots,d_N)=(\xi_1^{c_1}d_1,\ldots,\xi_N^{c_N}d_N),$$
for any $(c_1,\ldots,c_N)\in\Z^N$ and $(d_1,\ldots,d_N)\in\big(F^\times\big)^N$.  Form the semidirect product $G^r\rtimes S_r$ by letting the symmetric group $S_r$ act on $G^r$ (on the left) by permuting the factors of $G^r$:
$$^\pi (\rho_1,\ldots,\rho_r)=(\rho_{\pi(1)},\ldots,\rho_{\pi(r)}),$$
for all $\pi\in S_r$ and $\rho_i\in G$.
  This semidirect product acts on the space of ordered $r$-tuples of points in the torus $\big(F^\times\big)^N$ by letting $G^r$ act diagonally and letting $S_r$ permute the points:
\begin{equation}\label{512}
^{\rho\pi}a=(\rho_1.a_{\pi(1)},\ldots\rho_r.a_{\pi(r)}),
\end{equation}
for all $\rho=(\rho_1,\ldots,\rho_r)\in G^r$, $\pi\in S_r$, and $r$-tuples $a=(a_1,\ldots,a_r)$ of points $a_i\in\big(F^\times\big)^N$.

The group $G^r\rtimes S_r$ also acts on the space of $r$-tuples $\gl$ of nonzero dominant integral weights.  For each $\rho=(\rho_1,\ldots,\rho_r)\in G^r$ write $\rho_i=(\gs_1^{\rho_{i1}},\ldots,\gs_N^{\rho_{iN}})$ for some nonnegative integers $\rho_{ij}$.  Let the $\rho_i$ act on $\fg$ by
$$\rho_i(x)=\gs_1^{\rho_{i1}}\cdots\gs_N^{\rho_{iN}}x,$$
for all $x\in\fg$, and on the weights $\gl_i$ by
$$\rho_i(\gl_i)=\gl_i\circ\gamma(\rho_i^{-1}),$$
where $\gamma(\rho_i^{-1})$ is the outer part of the automorphism $\rho_i^{-1}:\ \fg\rightarrow\fg$.  Then $G^r\rtimes S_r$ acts on each $\gl=(\gl_1,\ldots,\gl_r)$ by
\begin{equation}\label{513}
^{\rho\pi}\gl=(\gl_{\pi(1)}\circ\gamma(\rho_1^{-1}),\ldots,\gl_{\pi(r)}\circ\gamma(\rho_r^{-1})).
\end{equation}

Combining (\ref{512}) and (\ref{513}) gives an action of $G^r\rtimes S_r$ on the set of pairs $(\gl,a)$, where $\gl$ is an $r$-tuple of nonzero dominant integral weights $\gl_i$ and $a$ is an $r$-tuple of points $a_i\in\big(F^\times\big)^N$ with $m(a_i)\neq m(a_j)$ whenever $i\neq j$.  Namely, let
\begin{equation}
^{\rho\pi}(\gl,a)=\left(^{\rho\pi}\gl,^{\rho\pi}a\right).
\end{equation}

In terms of this action, the isomorphism classes of the finite-dimensional simple $\cL$-modules are labelled by the orbits of the groups $G^r\rtimes S_r$.

\begin{corollary}\label{orbitcorollary}
Let $\gl=(\gl_1,\ldots,\gl_r)$ and $\mu=(\mu_1,\ldots,\mu_s)$ be sequences of nonzero dominant integral weights with respect to $\Delta$.  Suppose $a=(a_1,\ldots,a_r)$ and $b=(b_1,\ldots,b_s)$ are sequences of points in $\left(F^\times\right)^N$ with $m(a_i)\neq m(a_j)$ and $m(b_i)\neq m(b_j)$ whenever $i\neq j$.  

Then $V(\gl,a)$ and $V(\mu,b)$ are isomorphic if and only if $r=s$ and $(\gl,a)=^{\rho\pi}(\mu,b)$ for some $\rho\in G^r$ and $\pi\in S_r$.  
\end{corollary}
\proof Note that $m(a_i)=m(b_{\pi(i)})$ if and only if the coordinates $a_{ij}$ of $a_i=(a_{i1},\ldots,a_{iN})$ differ from the coordinates $b_{\pi(i)j}$ of $b_{\pi(i)}$ by an $m_j$th root of unity.  Since $\xi_j$ is a primitive $m_j$th root of unity, this happens if and only if there are integers $\rho_{ij}$ so that $a_{ij}=\xi_j^{\rho_{ij}}b_{\pi(i)j}$.  In terms of group actions, this is precisely the existence of $\rho_i=(\gs_1^{\rho_{i1}},\ldots,\gs_N^{\rho_{iN}})\in G$ with $a_i=\rho_i.b_{\pi(i)}$.  In other words, $m(a_i)=m(b_{\pi(i)})$ for all $i$ if and only if  $a=^{\rho\pi}b$ for some $\rho\in G^r$ and $\pi\in S_r$.

Since $\xi_j^{\rho_{ij}}=\frac{a_{ij}}{b_{\pi(i)j}}$, we see that 
\begin{align*}
\rho_i^{-1}(x)&=\gs_1^{-\rho_{i1}}\cdots\gs_N^{-\rho_{iN}}x\\
&=\xi_1^{-\rho_{i1}k_1}\cdots\xi_N^{-\rho_{iN}k_N}x\\
&=\frac{b_{\pi(i)}^k}{a_i^k}x,
\end{align*}
for all $k\in\Z^N$ and $x\in\fg_{\ol{k}}$.  Therefore, the automorphism $\omega_i$ of Theorem \ref{isothm} is equal to $\rho_i^{-1}$, and $\gl=^{\rho\pi}\mu$ is equivalent to the condition that $\gl_i=\mu_{\pi(i)}\circ\gamma_i$ for every $i$.\qed

\bigskip

For any diagram automorphism $\gs_1$, the finite-dimensional simple modules for the twisted (single) loop algebra $\cL(\fg;\gs_1)$ were classified in \cite{ChFoSe}.  Recently, Senesi, Neher, and Savage \cite{NeSaSe} have reinterpreted this work in terms of finitely supported $\gs_1$-equivariant maps $F^\times\rightarrow P_+$, where $P_+$ is the set of nonzero dominant integral weights of $\fg$ with respect to a fixed Cartan subalgebra and base of simple roots.  Theorem \ref{isothm} and Corollary \ref{orbitcorollary} can be used to extend this perspective to the multiloop setting.
 
Let $\gl=(\gl_1,\ldots,\gl_r)$ and $a=(a_1,\ldots,a_r)$ be as in Theorem \ref{isothm}.  Each evaluation module $V_{\gl_i}(a_i)$ corresponds to a map
\begin{align}
\chi_{\gl_i,a_i}:\ (F^\times)^N&\rightarrow P_+\\
x&\mapsto \gd_{x,a_i}\gl_i.
\end{align}
The isomorphism class $[\gl,a]$ of the tensor product $V(\gl,a)$ can then be identified with the sum of all the characters $\chi_{\eta_0,c_0}$ for which the pair $(\eta_0,c_0)=(\mu_1,b_1)$ for some $\mu=(\mu_1,\ldots,\mu_r)$ and $b=(b_1,\ldots,b_r)$ with $(\mu,b)$ in the $G^r\rtimes S_r$-orbit of $(\gl,a)$.  That is, we let
\begin{equation}
\chi_{[\gl,a]}=\sum_{g\in G}\sum_{i=1}^r\chi_{\gl_i\circ\gamma(g^{-1}),g.a_i}.
\end{equation}
This process associates a finitely supported $G$-equivariant map 
\begin{equation}
\chi_{[\gl,a]}:\ (F^\times)^N\rightarrow P_+
\end{equation}
with each isomorphism class of finite-dimensional simple $\cL(\fg;\gs_1,\ldots,\gs_N)$-modules.  From Corollary \ref{orbitcorollary} and the construction of $\chi_{[\gl,a]}$, it is easy to see that distinct isomorphism classes get sent to distinct functions.  

Conversely, any finitely supported $G$-equivariant map $f:\ (F^\times)^N\rightarrow P_+$ corresponds to an isomorphism class $[\gl,a]$ of finite-dimensional simple $\cL$-modules as follows.  By $G$-equivariance, the support $\hbox{supp}\,f$ of $f$ decomposes into a disjoint union of $G$-orbits.  Choose representatives $a_1,\ldots,a_r\in (F^\times)^N$ to label each $G$-orbit in $\hbox{supp}\,f$.  Since the $G$-orbits are disjoint, $m(a_i)\neq m(a_j)$ whenever $i\neq j$, and by definition of $f$,  $\gl:=(f(a_1),\ldots,f(a_r))$ is an $r$-tuple of nonzero dominant integral weights.  Then by Theorem \ref{evsaresimple}, $V(\gl,a)$ is a finite-dimensional simple $\cL$-module, and by Corollary \ref{orbitcorollary}, the isomorphism class $[f]:=[\gl,a]$ of this module is independent of the choice of orbit representatives $a_1,\ldots,a_r$.  It is now straightforward to verify that $\chi_{[f]}=f$ 
\ for all finitely supported $G$-equivariant maps $f:\ (F^\times)^N\rightarrow P_+$.

\begin{corollary}\label{equivariantmapcorollary}
The isomorphism classes of the finite-dimensional simple $\cL$-modules are in bijection with the finitely supported $G$-equivariant maps $(F^\times)^N\rightarrow P_+$.
\end{corollary}
\qed


\begin{thebibliography}{}


\bibitem{ABFP1}
B. Allison, S. Berman, J. Faulkner, and A. Pianzola, Realization of graded-simple algebras as loop algebras, Forum Math. {\bf 20} (2008), 395--432.

\bibitem{ABFP2}
B. Allison, S. Berman, J. Faulkner, and A. Pianzola, Multiloop realization of extended affine Lie algebras and Lie tori, preprint arXiv: math.RA/0709.0975v2.

\bibitem{AtiMac}
M. Atiyah and I. MacDonald, {\em Introduction to Commutative Algebra}, Addison-Wesley, London, 1969.

\bibitem{batra}
P. Batra, Representations of twisted multi-loop Lie algebras, J. Algebra {\bf 272} (2004), 404--416.




\bibitem{bourbaki}
N. Bourbaki, {\em \'El\'ements de math\'ematique: Alg\`ebre}, Chap\^itre 8, Hermann, Paris, 1958. 

\bibitem{bourbakiLie}
N. Bourbaki, {\em \'El\'ements de math\'ematique: Groupes et Alg\`ebres de Lie}, Chap\^itre 1, Hermann, Paris, 1960.

\bibitem{chari86}
V. Chari, Integrable representations of affine Lie-algebras, Invent. Math. {\bf 85} (1986), 317--335.

\bibitem{cp86}
V. Chari and A. Pressley, New unitary representations of loop groups, Math. Ann. {\bf 275} (1986), 87--104.


\bibitem{cp87}
V. Chari and A. Pressley, A new family of irreducible, integrable modules for affine Lie algebras, Math. Ann. {\bf 277} (1987), 543--562.


\bibitem{cp88}
V. Chari and A. Pressley, Integrable representations of twisted affine Lie algebras, J. Algebra {\bf 113} (1988), 438--464.

\bibitem{ChFoSe}
V. Chari, G. Fourier, and P. Senesi, Weyl modules for the twisted loop algebras, J. Algebra {\bf 319} (2008), 5016--5038.

\bibitem{GiPi}
P. Gille and A. Pianzola, Galois cohomology and forms of algebras over Laurent polynomial rings,  Math. Ann. {\bf 338} (2007), 497--543.

\bibitem{jacobson}
N. Jacobson, {\em Lie Algebras}, Dover, New York, 1962.


\bibitem{kacLiebook}
V. Kac, {\em Infinite Dimensional Lie Algebras}, 3rd edition, Cambridge University Press, Cambridge, 1990.


\bibitem{neher}
E. Neher, Extended affine Lie algebras, C.R. Math. Acad. Sci. Soc. R. Can. {\bf 26} (2004), 90--96.

\bibitem{NeSaSe}
E. Neher, A. Savage, and P. Senesi, private communication.



\bibitem{rao93}
S.E. Rao, On representations of loop algebras, Comm. Algebra {\bf 21} (1993), 2131--2153.

\bibitem{rao01}
S.E. Rao, Classification of irreducible integrable modules for multi-loop algebras with finite-dimensional weight spaces, J. Algebra {\bf 246} (2001), 215--225.




\end{thebibliography}
\end{document}